\documentclass[10pt,leqno]{article}

\usepackage{amsmath,amssymb,amsthm,mathrsfs,dsfont,graphicx}

\usepackage[margin=3cm]{geometry} 

\usepackage{titlesec,hyperref}
\usepackage{esint}
\usepackage{natbib}
\usepackage{color}

\usepackage{fancyhdr}
\titleformat{\subsection}{\it}{\thesubsection.\enspace}{1pt}{}

\newtheorem{theo}{Theorem}[section]
\newtheorem{lemm}[theo]{Lemma}
\newtheorem{defi}[theo]{Definition}

\newtheorem{prop}[theo]{Proposition}
\newtheorem{rema}[theo]{Remark}
\numberwithin{equation}{section}

\allowdisplaybreaks 

\begin{document}
\title{H\"{o}lder continuous weak solutions of the 3D Boussinesq equation with thermal diffusion
\hspace{-4mm}
}

\author{Zipeng $\mbox{Chen}^1$ \footnote{Email: chenzp26@mail2.sysu.edu.cn},\quad
	 and\quad
	Zhaoyang $\mbox{Yin}^{1}$\footnote{E-mail: mcsyzy@mail.sysu.edu.cn}\\
    $^1\mbox{School}$ of Science,\\ Shenzhen Campus of Sun Yat-sen University, Shenzhen 518107, China
	}
        
\date{}
\maketitle
\hrule

\begin{abstract}
In this paper, we show the existence of H\"{o}lder continuous periodic weak solutions of the 3D Boussinesq equation with thermal diffusion, which apprroximate the Onsager's critical spatial regularity and satisfy the prescribed kinetic energy. More precisely, for any smooth $e(t):[0,T]\rightarrow \mathbb{R}_+$ and $\beta\in (0, \frac{1}{3})$, there exist $v\in C^{\beta}([0,T]\times {\mathbb{T} }^3)$ and $ \theta\in C_t^{1,\frac{\beta}{2}}C_x^{2,\beta}([0,T]\times {\mathbb{T} }^3)$ which solve (\ref{e:boussinesq equation}) in the sense of distribution and satisfy
\begin{align}
e(t)=\int_{{{\mathbb{T} }^3}}|v(t,x)|^2dx, \quad \forall t\in [0,T].\nonumber
\end{align} 
\end{abstract}
\noindent {\sl Keywords:}  Boussinesq equation, H\"{o}lder continuous periodic weak solutions, Convex integration, Onsager exponent

\vskip 0.2cm

\noindent {\sl AMS Subject Classification:} 35Q30, 76D03  \

\vspace*{10pt}


\section{Introduction and main result}

  In this paper, we consider the following 3D Boussinesq equation with thermal diffusion:
\begin{equation}\label{e:boussinesq equation}
\begin{cases}
\partial_tv+ \text{div}(v\otimes v)+\nabla p=\theta e_3, \quad\quad \\
\text{div}\,v=0,\\
\partial_t\theta+v\cdot \nabla\theta-\Delta \theta=0, \quad\forall(t,x)\in [0,T]\times\mathbb{T}^3
\end{cases}
\end{equation}
where $T>0$, $\mathbb{T}^3$ is the 3-dimensional torus and $e_3=(0,0,1)$. Here, $v,p,\theta$ represent velocity, pressure, and temperature, respectively. The Boussinesq equation was introduced to model large-scale atmospheric and oceanic flows that are responsible for cold fronts and the jet stream (see, \cite{Majda}). 

The study of weak solutions in fluid dynamics, including those that fail to conserve kinetic energy, has been popular in recent years. One of the famous problems is Onsager's conjecture, which states that the incompressible Euler equation admits H\"{o}lder continuous weak solutions that dissipate the kinetic energy. More precisely, Lars Onsager conjectured:
   
\noindent(1) For any $\alpha >\frac{1}{3}$, every $C_{t,x}^{\alpha}$ weak solution conserves energy.

\noindent(2) For any $\alpha <\frac{1}{3}$, there exist dissipative solutions with $C_{t,x}^{\alpha}$ regularity.

  Part (1) of this conjecture was proved by Constantin et al.\cite{cons}. Duchon and Robert \cite{Duchon} and Cheskidov et al.\cite{Che} gave a proof of this part with weaker assumptions on the solution. After some improvement on Part (2) \cite{1,2,3,4,5,6,7,8,9}, it was finally proved by Isett \cite{Isett3} (i.e., reaching a critical regularity $\frac{1}{3}$). Later, Buckmaster et al. gave another proof in \cite{B1}. Furthermore, the idea and technique are utilized to construct dissipative weak solutions or prove nonuniqueness for other equations (see, e.g., \cite{ns不唯一}, \cite{SQG}, \cite{activescalar}, and \cite{quasi-geostrophic}).

  For the 3D inviscid Boussinesq equation, Tao and Zhang \cite{T3} showed the existence of $C^{\alpha}$ periodic weak solutions with the prescribed kinetic energy where $\alpha\in(0,\frac{1}{5})$.  Miao et al.\cite{Ye} proved the H\"{o}lder threshold regularity exponent for $L^p$-norm conservation of temperature of this system is $\frac{1}{3}$. Xu and Tan  \cite{xusaigu} demonstrated the existence of $C^{\alpha}$ periodic weak solutions with the prescribed kinetic energy where $\alpha\in(0,\frac{1}{3})$. In the two-dimensional case, Luo et al.\cite{L2,L1} constructed the H\"{o}lder continuous dissipative weak solutions of the Boussinesq equation with fractional dissipation in velocity and thermal diffusion.

  Motivated by the work of \cite{B1} and \cite{L1}, we construct the H\"{o}lder continuous dissipative weak solutions of (\ref{e:boussinesq equation}) by combining the convex integration method and the energy method. More precisely, the velocity is constructed by the Mikado flow, as usual. In order to overcome the difficulty of interactions between velocity and temperature and the impact of thermal diffusion, the temperature is derived by solving the transport-diffusion equation. Unlike the convex integration method of \cite{L1}, which is based on the Beltrami flow and the H\"{o}lder exponent drops to $\frac{1}{10}$, we construct the $C^{\alpha}$ periodic weak solutions with the prescribed kinetic energy for any $\alpha\in(0,\frac{1}{3})$. In fact, the influence of the temperature effect on the velocity appears to be not strong enough to change the critical H\"{o}lder exponent $\frac{1}{3}$ in 3-dimensional case. Our main theorem is as follows.
   
   \begin{theo}\label{thm1}
   	Assume that $e(t): [0,T]\rightarrow \mathbb{R}_+$ is a strictly positive smooth function. Let $\theta^0(x_3)$ a smooth function only depended on $x_3$ with zero mean and $\alpha \in (0, \frac{1}{3})$. Then there exist $v \in C^{\alpha}([0,T]\times {\mathbb{T} }^3)$, $\theta \in C^{1,\frac{\alpha}{2}}_tC^{2,\alpha}_x$ such that $(v,\theta)$ satisfy (\ref{e:boussinesq equation}) in the sense of distribution and 
    \begin{gather}
   e(t)=\int_{{{\mathbb{T} }^3}}|v(t,x)|^2dx,\\ \frac12\|\theta(t,\cdot)\|_{L^2}^2+\int_0^t\|\nabla \theta(s,\cdot)\|^2_{L^2}ds=\frac12\|\theta^0(\cdot)\|_{L^2}^2.
    \end{gather} 
   \end{theo}

\begin{rema}
    In Theorem \ref{thm1}, when $\theta^0=0$,  the velocity we construct is indeed the solution of the Euler equation, which has been studied in \cite{B1}. Nevertheless, we are able to select any non-zero $\theta^0$ in a special class to construct a non-zero $\theta$. We denote two quantity:\begin{gather*}
        E(t)\triangleq\frac12\|v(t,\cdot)\|_{L^2}^2-\int_0^t\int_{{{\mathbb{T} }^3}}\theta_3v_3\, dx,\\
    M(t)\triangleq\frac12\|\theta(t,\cdot)\|_{L^2}^2+\int_0^t\|\nabla \theta(s,\cdot)\|^2_{L^2}ds.
    \end{gather*}
    $E(t)$ and $M(t)$ are constant if $(v,\theta)$ is a smooth solution of (\ref{e:boussinesq equation}). More precisely, we have an Onsager-type theorem for (\ref{e:boussinesq equation}) which extends the Onsager theorem of the Euler equation: \begin{itemize}
        \item For any $\alpha>\frac{1}{3}$, if $u\in C^{\alpha}([0,T]\times\mathbb{T}^3)$, then $E(t)$ and $M(t)$ are constant.
        \item For any $\alpha<\frac{1}{3}$, there exist $u\in C^{\alpha}([0,T]\times\mathbb{T}^3)$ and $\theta\in L^{\infty}(0,T;L^2(\mathbb{T}^3))\cap L^2(0,T;\dot{H^1}(\mathbb{T}^3))$ such that $E(t)$ is not constant and $M(t)$ remains constant.
        \end{itemize}
    In fact, the first result can be proved by using the method in \cite{cons} and the second is due to Theorem \ref{thm1}.

\end{rema}

\section{The proof of the main result}
Given an initial datum $\theta^0(x_3)$ with zero mean and for any $q\in\mathbb{N}$, we will construct a smooth solution $(v_q,\theta_q,\mathring{R_q})$ that satisfies the following Boussinesq-Reynolds equation on $[0,T]\times {\mathbb{T} }^3$:
    \begin{equation}\label{e:Boussinesq-Reynold}
    \begin{cases}
    \partial_tv_q+\text{div}(v_q\otimes v_q)+\nabla p_q=\theta_q e_3+\text{div}\mathring{R_q}, \quad\quad \\
    \text{div}\,v_q=0,\\
    \partial_t\theta_q+v_q\cdot \nabla\theta_q-\Delta \theta_q=0, \\ 
    \theta_q(0,x)=\theta^0(x_3),\quad\quad 
    \end{cases}
    \end{equation}
where $\mathring{R_q}$ is a symmetric matrix, moreover we add the constraints that
\begin{align}
    tr(\mathring{R}_q)=0\label{tr0}
\end{align}
and that
\begin{align}
    \int_{\mathbb{T}^3}p_q(t,x)dx=0\label{intp=0}.
\end{align}

For $q\in \mathbb{N}$, we define two parameters $\lambda_q$ and $\delta_q$ to measure the size of the corresponding solution:
\begin{align}
    \lambda_q=2\pi\lceil a^{(b^q)}\rceil, \delta_q=\lambda_q^{-2\beta}
\nonumber
\end{align}
where $\lceil x\rceil$ denotes the smallest integer $n\geq x$. In the proof, we will choose $a\gg 1$, $b>1$ is nearly equal to 1 and $0<\beta<\frac{1}{3}$. 

The following proposition serves as the basis for the proof of the Theorem \ref{thm1}.
\begin{prop}\label{Prop1}
   	Assume $0<\beta<\frac{1}{3}$ and $1<b<\frac{\beta+\sqrt{4\beta-3\beta^2}}{4\beta}$. Let $e(t)$ and $\theta^0(x_3)$ be as in Theorem \ref{thm1}. Then there exists $M$ and $C_0$ depending on $e(t)$, $\{C(N)\}_{N\geq2}$ depending on $M$, $\alpha$ depending on $\beta$ and $b$, and $a$ depending on $\beta,b,\alpha$, $M$ and $C_0$ such that the following holds: there exists a sequence of functions $(v_q,p_q,\theta_q,\mathring{R_q})\in C^\infty([0,T]\times \mathbb{T}^3)$ starting from $(v_0,p_0,\theta_0,\mathring{R_0})$, satisfying the (\ref{e:Boussinesq-Reynold})-(\ref{intp=0}) and the following estimates: 
    \begin{itemize}
        \item
    the Reynolds stress error $\{\mathring{R}_q\}$ satisfies
    \begin{gather}
    \|\mathring{R}_q\|_0\leq \delta_{q+1}\lambda_q^{-3\alpha}\label{p1},
    \end{gather}
    \item the velocity $\{v_q\}$ satisfies
    \begin{gather}
    \|v_q\|_0\leq C_0-\delta^{\frac{1}{2}}_q\label{p2},\\
    \|v_q\|_1\leq M\delta_q^{\frac{1}{2}}\lambda_q\label{p3},\\
    \|v_q\|_N\leq C(N)\delta_q^{\frac{1}{2}}\lambda_q^N, \quad\forall N\geq2\label{p4},\\
    \delta_{q+1}\lambda_q^{-\alpha}\leq e(t)-\int_{{\mathbb{T}}^3}|v_q(t,x)|^2dx\leq \delta_{q+1},\quad\forall t\in[0,T],\label{p5}\\
    \|v_{q+1}-v_{q}\|_0+\lambda_{q+1}^{-1}\|v_{q+1}-v_{q}\|_1\leq M\delta_{q+1}^{\frac{1}{2}}\label{p6},
    \end{gather}
    \item the temperature $\{\theta_q\}$ satisfies
    \begin{gather}
    \|(\theta_{q+1}-\theta_{q})(t,\cdot)\|^2_{L^2}+\int_0^t\|\nabla(\theta_{q+1}-\theta_{q})(s,\cdot)\|^2_{L^2}ds\leq C\delta_{q+1}^{\frac{1}{2}},\quad\forall t\in[0,T],\label{p7}\\
    \frac{1}{2}\|\theta_{q}(t,\cdot)\|^2_{L^2}+\int_0^t\|\nabla\theta_{q}(s,\cdot)\|^2_{L^2}ds=\frac{1}{2}\|\theta^0(\cdot)\|^2_{L^2},\quad\forall t\in[0,T].\label{p8}
    \end{gather}
    \end{itemize}
\end{prop}
The H\"{o}lder norms employed above are defined in Appendix \ref{appendix A}, which we consider only the spatial regularity. The proof of Proposition \ref{Prop1} will occupy most of the paper and will be presented later in this paper. This proposition immediately proves Theorem \ref{thm1} and the proof is similar to \cite{B1}, so we first give a simple sketch here.

\begin{proof}[Proof of Theorem\ref{thm1}]
    By (\ref{p1}), (\ref{p6}) and (\ref{p7}), we notice that $\{v_q\}$ converges uniformly to a continuous function $v$, $\{\mathring{R_q}\}$ converges uniformly to 0 and $\{\theta_q\}$ converges to a function $\theta$ in $L^{\infty}(0,T;L^2(\mathbb{T}^3))\cap L^2(0,T;\dot{H^1}(\mathbb{T}^3))$. Moreover, since we have
    \begin{align}
        \Delta p_q=\text{div}\text{div}(-v_q\otimes v_q+\mathring{R}_q)+\text{div}(\theta_qe_3)
    \end{align}
    and (\ref{intp=0}), we get $p_q$ also converges to some pressure $p$ in $L^{\infty}(0,T;L^2(\mathbb{T}^3))$. Passing to the limit in (\ref{e:Boussinesq-Reynold}), we deduce that $(v,p,\theta)$ satisfy (\ref{e:boussinesq equation}) in the sense of distribution.

    Using (\ref{p7}) and choosing $0<\beta'<\beta$, we infer that
    \begin{align*}
        \sum^{\infty}_{q=0}\|v_{q+1}-v_q\|_{\beta'}&\lesssim \sum^{\infty}_{q=0}\|v_{q+1}-v_q\|_0^{1-\beta'}\|v_{q+1}-v_q\|_1^{\beta'}\\
        &\lesssim\sum^{\infty}_{q=0}\delta_{q+1}^{\frac{1-\beta'}{2}}(\delta_{q+1}^{\frac{1}{2}}\lambda_{q+1})^{\beta'}\\
        &\lesssim\sum^{\infty}_{q=0}\lambda_{q+1}^{\beta'-\beta}.
    \end{align*}
    Here and throughout the paper, $x\lesssim y$ denotes $x\leq Cy$ for a constant $C>0$ that is independent of $a,b$ and $q$, but may change from line to line. Hence, we obtain $v\in C^0_tC^{\beta'}_x$ for all $\beta'<\beta$. 

    Since $\theta$ satisfies
    \begin{align*}
        \begin{cases}
            \partial_t\theta+v\cdot \nabla\theta-\Delta \theta=0,\\
            \theta(0,x)=\theta^0(x_3),
        \end{cases}
    \end{align*}
    we deduce that $\|\theta\|_0\leq||\theta^0\|_0$ by the maximum principle.

    Next we give a proof of recovering the time regularity of $v$. Let $\tilde{v}_q=v*\varphi_{2^{-q}}$, where $\varphi$ is a smooth standard mollifier in space and $\varphi_l(x)=l^{-3}\varphi(xl^{-1})$. Using standard mollification estimates, we have
    \begin{align}
        \|\tilde{v}_q-v\|_0\lesssim\|v\|_{\beta'}2^{-q\beta'}\lesssim2^{-q\beta'}\label{t1}.
    \end{align}
    Moreover, $\tilde{v}_q$ satisfies the equation
    \begin{align*}
         \partial_t\tilde{v}_q+\text{div}(v\otimes v)*\varphi_{2^{-q}}+\nabla p*\varphi_{2^{-q}}=\theta e_3*\varphi_{2^{-q}}.
    \end{align*}
    Next, since
    \begin{align*}
        \Delta p*\varphi_{2^{-q}}=-\text{div}\text{div}(v\otimes v)*\varphi_{2^{-q`}}+\text{div}\,(\theta e_3)*\varphi_{2^{-q}},
    \end{align*}
    for any fixed $\epsilon>0$ the Schauder's estimates yield
    \begin{align*}
        \|\nabla p*\varphi_{2^{-q}}\|_{\epsilon}\lesssim\|v\otimes v\|_{\beta'}2^{q(1+\epsilon-\beta')}+\|\theta||_0 2^{q\epsilon}\lesssim 2^{q(1+\epsilon-\beta')} .
    \end{align*}
    Moreover, 
    \begin{align*}
       \|\text{div}(v\otimes v)*\varphi_{2^{-q}}\|_0\lesssim \|v\otimes v\|_{\beta'}2^{q(1-\beta')}\lesssim 2^{q(1-\beta')}.
    \end{align*}
    Thus, we deduce that
    \begin{align}
        \|\partial_t\tilde{v}_q\|_0\lesssim 2^{q(1+\epsilon-\beta')}\label{t2}.
    \end{align}
    For any $\beta''<\beta'$, choosing $\epsilon>0$ sufficiently small such that $\beta'-(1+\epsilon)\beta''\geq \epsilon$, we obtain from (\ref{t1}) and (\ref{t2}) that 
    \begin{align*}
        \|\tilde{v}_{q+1}-\tilde{v}_{q}\|_{C^{\beta''}_tC^0_x}&\lesssim(\|\tilde{v}_{q+1}-v\|_0+\|\tilde{v}_q-v\|_0)^{1-\beta''}(  \|\partial_t\tilde{v}_{q+1}\|_0+\|\partial_t\tilde{v}_q\|_0)^{\beta''}\\
        &\lesssim2^{-q\beta'(1-\beta'')+q\beta''(1+\epsilon-\beta')}\\
        &\lesssim2^{-q(\beta'-(1+\epsilon)\beta'')}\\
        &\lesssim2^{-q\epsilon}.
    \end{align*}
    Thus, the series
    \begin{align*}
        v=\tilde{v}_0+\sum_{q\geq0}(\tilde{v}_{q+1}-\tilde{v}_{q})
    \end{align*}
    converges in $C^{\beta''}_tC^0_x$. Combined with $v\in C^0_tC^{\beta'}_x$, we get $v\in C^{\beta''}([0,T]\times\mathbb{T}^3)$ with $\beta''<\beta'<\beta<\frac{1}{3}$ . Moreover, by the Schauder estimate of linear parabolic equation we have that $\theta\in C_t^{1,\frac{\beta''}{2}}C_x^{2,\beta''}$ with $\beta''<\beta<\frac{1}{3}$.

    Finally, let $q\to \infty$, from (\ref{p5}) and (\ref{p8}) we have
    \begin{gather*}
        e(t)=\int_{\mathbb{T}^3}|v(t,x)|^2dx,\quad \forall t\in [0,T],\\
        \frac{1}{2}\|\theta(t,\cdot)\|^2_{L^2}+\int_0^t\|\nabla\theta(s,\cdot)\|^2_{L^2}ds=\frac{1}{2}\|\theta^0(\cdot)\|^2_{L^2},\quad \forall t\in[0,T],
    \end{gather*}
    which completes the proof of Theorem \ref{thm1}.
    \end{proof}

The rest of the paper focuses on the proof of Proposition \ref{Prop1}. As in \cite{B1}, we construct ${v_q}$ using the inductive procedure and the convex integration scheme. Roughly speaking, there are three steps to construct $v_{q+1}$ from $v_q$: mollification, gluing, and perturbation. After constructing the new velocity $v_{q+1}$, we construct the new temperature $\theta_{q+1}$ by directly solving the transport-diffusion equation:
\begin{align*}
    \begin{cases}
            \partial_t\theta_{q+1}+v_{q+1}\cdot \nabla\theta_{q+1}-\Delta \theta_{q+1}=0,\\
            \theta_{q+1}(0,x)=\theta^0(x_3),
        \end{cases}
\end{align*}
Finally, we construct $\mathring{R}_{q+1}$ such that $(v_{q+1},p_{q+1},\mathring{R}_{q+1},\theta_{q+1})$ solve the equation (\ref{e:Boussinesq-Reynold}) and satisfy (\ref{p1})-(\ref{p8}) with $q$ replaced by $q+1$.

\section{The choice of starting cases}\label{start}
Firstly, we consider the construction of $(v_0, p_0,\mathring{R}_0,\theta_0)$. Unlike the Euler equation, the transformation 
\begin{align*}
    v(t,x)\mapsto\Gamma v(\Gamma t,x)
\end{align*}
is not applicable to the Boussinesq equation with thermal diffusion (\ref{e:boussinesq equation}), so instead of setting $(v_0, p_0,\mathring{R}_0)=(0,0,0)$ and making further assumption on $e(t)$, we choose the starting vector as follows: 

\begin{gather*}
v_0=
    \begin{pmatrix}
        \sqrt{\frac{2e(t)-\delta_1-\delta_1\lambda_0^{-\alpha}{}}{8\pi^3}}sin(\lceil\delta_0^{\frac{1}{2}}\lambda_0\rceil x_2)\\0\\0
    \end{pmatrix}\\
\mathring{R}_0=
    \begin{pmatrix}
        0&-\frac{e'(t)}{\sqrt{8\pi^3(2e(t)-\delta_1-\delta_1\lambda_0^{-\alpha})}}\frac{cos(\lceil\delta_0^{\frac{1}{2}}\lambda_0\rceil x_2)}{\lceil\delta_0^{\frac{1}{2}}\lambda_0\rceil}&0\\
        -\frac{e'(t)}{\sqrt{8\pi^3(2e(t)-\delta_1-\delta_1\lambda_0^{-\alpha})}}\frac{cos(\lceil\delta_0^{\frac{1}{2}}\lambda_0\rceil x_2)}{[\delta_0^{\frac{1}{2}}\lambda_0]}&0&0\\
        0&0&0
    \end{pmatrix}\\
    p_0=\int_0^{x_3}e^{t\Delta}\theta^0(s)ds-f(t), \quad\theta_0=e^{t\Delta}\theta^0(x_3).
\end{gather*}
where $f(t)$ is a function such that $\int_{\mathbb{T}^3}p_0=0$.

We note 
\begin{align}
    M_1=\sup_t \{|e(t)|+|e'(t)|\},  \quad m_1=\inf_t e(t)\label{M1}.
\end{align}
Since $b<\frac{\beta+\sqrt{4\beta-3\beta^2}}{4\beta}<\frac{1-\beta}{2\beta}$, $(v_0, p_0,\mathring{R}_0, \theta_0)$ satisfies
\begin{gather}
    \|v_0\|_0\leq\sqrt{\frac{M_1}{4\pi^3}}\leq \sqrt{\frac{M_1}{4\pi^3}}+1-\delta_0^{\frac{1}{2}},\\
    \|v_0\|_N\leq\sqrt{\frac{M_1}{4\pi^3}}\delta_0^{\frac{1}{2}}\lambda_0^N, \quad\forall N\geq1,\label{v0N}\\
    \|\mathring{R}_0\|_0\leq\frac{M_1}{\sqrt{8\pi^3m_1}\delta_0^{\frac{1}{2}}\lambda_0}\leq\delta_1\lambda_0^{-3\alpha}\label{R0},\\
    \delta_{1}\lambda_0^{-\alpha}\leq e(t)-\int_{{\mathbb{T}}^3}|v_0(t,x)|^2dx\leq \delta_{1},\quad\forall t\in[0,T],\\
    \frac{1}{2}\|\theta_{0}(t,\cdot)\|^2_{L^2}+\int_0^t\|\nabla\theta_{0}(s,\cdot)\|^2_{L^2}ds=\frac{1}{2}\|\theta^0(\cdot)\|^2_{L^2},\quad\forall t\in[0,T].
\end{gather}
Here we use $\delta_1-\delta_1\lambda_0^{-\alpha}\leq m_1$ and $(\delta_0^{\frac{1}{2}}\lambda_0)^{-1}\lesssim\delta_1\lambda_0^{-3\alpha}$ if a is sufficiently large and $\alpha$ is sufficiently small.
Thus we set the constant $C_0$ in Proposition \ref{Prop1}:
\begin{align*}
    C_0=\sqrt{\frac{M_1}{4\pi^3}}+1.
\end{align*}
Moreover, $(v_0, p_0,\mathring{R}_0, \theta_0)$ obviously satisfy (\ref{e:Boussinesq-Reynold}), (\ref{tr0}) and (\ref{intp=0}) with $q=0$ . 

In the following, we will show the construction of $(v_{q+1}, p_{q+1},\mathring{R}_{q+1},\theta_{q+1})$ from $(v_q, p_q,\mathring{R}_q,\theta_q)$. Assuming $(v_q, p_q,\mathring{R}_q,\theta_q)$ satisfy (\ref{e:Boussinesq-Reynold}), (\ref{p1}), (\ref{p2}),(\ref{p3}) and (\ref{p5}), we first construct $v_{q+1}$ by using convex integration schemes.
\begin{rema}
    The constant M in (\ref{p3}) and (\ref{p6}) will be chosen in Section \ref{Perturbation}. Moreover, $\{C(N)\}_{N\geq2}$ in (\ref{p4}) will be unnecessary in the proof of constructing $v_{q+1}$, but it will be calculated directly in the Section \ref{Perturbation}. In fact, in Section \ref{温度构造}, assumptions (\ref{p3}) and (\ref{p4}) are devoted to deduce the estimates of $\theta_q$ in the Sobolev space.
\end{rema}

\section{Mollification}

Let $\phi$ be a standard mollifier in space, we set
\begin{align}
    l=\frac{\delta_{q+1}^{\frac{1}{2}}}{\delta_q^{\frac{1}{2}}\lambda_q^{1+\frac{3\alpha}{2}}}\label{l}.
\end{align}
Choosing $a$ sufficiently large and $\alpha$ sufficiently small, we get 
\begin{align}
    \lambda_q^{-\frac{3}{2}}\leq l \leq\lambda_q^{-1}.\label{l的范围}
\end{align}
We define 
\begin{align*}
    v_l=v_q*\phi_l, \quad\mathring{R}_l=\mathring{R}_q*\phi_l-(v_q \mathring{\otimes} v_q)*\phi_l+v_l \mathring{\otimes} v_l, \quad\theta_l=\theta_q*\phi_l,
\end{align*}
where $f\mathring{\otimes}g$ denotes the traceless part of $f\otimes g$. Since (\ref{e:Boussinesq-Reynold}), we have
\begin{align}
    \begin{cases}
\partial_tv_l+\text{div}(v_l\otimes v_l)+\nabla p_l=\theta_l e_3+\text{div}\mathring{R}_l,\\
\text{div}\,v_l=0,\label{e:vl}\\
    \end{cases}
\end{align}
for some suitable $p_l$.
Using (\ref{l的范围}), standard mollification estimates and (\ref{commu}), we can easily obtain the following proposition. The proof can be found in \cite[Proposition 2.2]{B1}.
\begin{prop}\label{propvl-vq}
    \begin{align}
        \|v_l-v_q\|_0&\lesssim\delta_{q+1}^{\frac{1}{2}}\lambda_q^{-\alpha},\label{vl-vq0}\\
        \|v_l\|_{N+1}&\lesssim\delta_q^{\frac{1}{2}}\lambda_ql^{-N}, \quad\forall N\geq0,\label{vl}\\
        \|\mathring{R}_l\|_{N+\alpha}&\lesssim\delta_{q+1}l^{-N+\alpha}, \quad\forall N\geq0,\label{Rl}\\
        \|\theta_l\|_N&\lesssim l^{-N}, \quad\forall N\geq0\label{thetal},\\
        |\int_{\mathbb{T}^3}|v_q|^2-|v_l|^2dx|&\lesssim\delta_{q+1}l^{\alpha}.
    \end{align}
\end{prop}

\section{Gluing}
\subsection{Estimates for classical Exact Solutions}
We introduce a temporal parameter:
\begin{align}
    \tau_q=\frac{l^{2\alpha}}{\delta_q^{\frac{1}{2}}\lambda_q}
\end{align}
and we define $t_i=i\tau_q$ for each $i\in \mathbb{N}$ such that $i\tau_q\leq T$. We consider the following Euler equation with a given external force:
\begin{align}
    \begin{cases}\label{Euler with F}
        \partial_tv_i+\text{div}(v_i\otimes v_i)+\nabla p_i=\theta_l e_3,\\
        \text{div}\, v_i=0,\\
        v_i(t_i,\cdot)=v_l(t_i,\cdot).
    \end{cases}
\end{align}
The proof of the existence of a unique solution of (\ref{Euler with F}) is standard (see, e.g., \citep[chap.7]{bcd}), so we omit the proof here. We focus on the proof of the following proposition.

\begin{prop}
    Let $a$ to be sufficiently large and $v_i$ satisfies (\ref{Euler with F}), for $t\in[t_i-\tau_q,t_i+\tau_q]$, we have
    \begin{align}
        \|v_i(t,\cdot)\|_{N+\alpha}\lesssim\delta_q^{\frac{1}{2}}\lambda_ql^{1-N-\alpha}\lesssim\tau_q^{-1}l^{1-N+\alpha},\quad\forall N\geq 1.\label{vi}
    \end{align}
\end{prop}

\begin{proof}
    Let $N\geq1$ and $\gamma$ be a multi-index with $|\gamma|=N$. Using $\Delta p_i=-tr(\nabla v_i\nabla v_i)+\text{div}(\theta_l e_3)$, (\ref{thetal}) and Schauder estimates, we get
    \begin{align*}
        \|\nabla \partial^{\gamma}p_i\|_\alpha\lesssim\|-tr(\nabla v_i\nabla v_i)+\text{div}(\theta_l e_3)\|_{N-1+\alpha}\lesssim\|v_i\|_{1+\alpha}\|v_i\|_{N+\alpha}+l^{-N-\alpha}.
    \end{align*}
    Using interpolation inequality in H\"{o}lder space,
    \begin{align*}
        \|[\partial^\gamma,v_i\cdot\nabla]v_i\|_\alpha\lesssim\|v_i\|_{1+\alpha}\|v_i\|_{N+\alpha}.
    \end{align*}
    Since
    \begin{align*}
        \partial_t\partial^{\gamma}v_i+v_i\cdot\nabla\partial^{\gamma}v_i+[\partial^\gamma,v_i\cdot\nabla]v_i+\nabla \partial^{\gamma}p_i=\partial^{\gamma}\theta_l e_3,
    \end{align*}
    thus we deduce
    \begin{align}
        \|(\partial_t+v_i\cdot\nabla)\partial^{\gamma}v_i\|_\alpha\lesssim\|v_i\|_{1+\alpha}\|v_i\|_{N+\alpha}+l^{-N-\alpha}\label{vi2}.
    \end{align}
    By applying (\ref{trans1}) and (\ref{vl}), for $N=1$, we obtain
    \begin{align*}
    \|v_i(t)\|_{1+\alpha}&\lesssim(\|v_l(t_i)\|_{1+\alpha}+\int_{t_i}^t\|v_i(s)\|_{1+\alpha}^2+l^{-1-\alpha}ds)exp(\int_{t_i}^t\|v_i(s)\|_{1+\alpha}ds)\\
        &\lesssim(\delta_q^{\frac{1}{2}}\lambda_ql^{-\alpha}+\tau_ql^{-1-\alpha}+\int_{t_i}^t\|v_i(s)\|_{1+\alpha}^2ds)exp(\int_{t_i}^t\|v_i(s)\|_{1+\alpha}ds).
    \end{align*}
    Hence, by basic connectivity arguments and the fact $\tau_ql^{-1-\alpha}\leq\delta_q^{\frac{1}{2}}\lambda_ql^{-\alpha}$ if a is sufficiently large and $\alpha$ is sufficiently small, we get
    \begin{align}
        \|v_i(t)\|_{1+\alpha}\lesssim\delta_q^{\frac{1}{2}}\lambda_ql^{-\alpha},\quad\forall |t-t_i|\leq\tau_q\label{v1}.
    \end{align}
    Finally, (\ref{vi}) follows as a consequence of (\ref{vl}), (\ref{vi2}),(\ref{v1}), (\ref{trans1}) and Gr\"{o}nwall's inequality.
\end{proof}

Due to (\ref{e:vl}) and (\ref{Euler with F}), we have
    \begin{equation}\label{e:vl-vi}
    \begin{cases}
    \partial_t(v_l-v_i)+v_l\cdot\nabla(v_l-v_i)=(v_i-v_l)\cdot\nabla v_i-\nabla(p_l-p_i)+\text{div}\mathring{R}_l,\\
    v_l-v_i|_{t=t_i}=0,
    \end{cases}
    \end{equation}
thus $v_i-v_l$ is small when $|t-t_i|$ is small. Moreover, from the identity
$v_i-v_{i+1}=(v_i-v_l)-(v_{i+1}-v_l)$, we know $v_i$ is also close to $v_{i+1}$. In fact, we have the following estimates:
\begin{prop}
    For $|t-t_i|\leq\tau_q$ and $N\geq0$, it holds that
    \begin{gather}
    \|v_i-v_l\|_{N+\alpha}\lesssim\tau_q\delta_{q+1}l^{-N-1+\alpha},\label{vi-vl}\\
    \|\nabla(p_l-p_i)\|_{N+\alpha}\lesssim\delta_{q+1}l^{-N-1+\alpha},\\
    \|D_{t,l}(v_i-v_l)\|_{N+\alpha}\lesssim\delta_{q+1}l^{-N-1+\alpha},
    \end{gather}
    where we define the transport derivative
    \begin{align*}
        D_{t,l}=\partial_t+v_l\cdot\nabla.
    \end{align*}
\end{prop}
The proof of the above proposition is based on (\ref{vl}), (\ref{Rl}) and (\ref{vi}) and (\ref{trans1}), which can be found in \cite[Proposition 3.3]{B1}. 

Next, we define the first-order potentials:
\begin{align*}
    z(v)\triangleq\mathcal{B}v=(-\Delta)^{-1}\text{curl}\,v,
\end{align*}
where $\mathcal{B}$ is called Biot-Savart operator and $v$ is a smooth vector function on $\mathbb{T}^3$. Moreover, we have
\begin{align*}
    \text{div} \, z(v)=0\quad and\quad \text{curl}\,z(v)=v-\fint_{\mathbb{T}^3}v\,dx.
\end{align*}
For any $i$ we denote 
\begin{align*}
    z_i=z(v_i)\quad \text{and}\quad \tilde{z}_i=z(v_i-v_l).
\end{align*}
Our goal is to obtain estimates for $\tilde{z}_i$ and $z_i-z_{i+1}$. From (\ref{vi-vl}), we have
\begin{align*}
    \|v_i-v_{i+1}\|_{N+\alpha}\lesssim\tau_q\delta_{q+1}l^{-N-1+\alpha},
\end{align*}
so we expect to obtain a factor $l$ since the characteristic length of $v_i-v_{i+1}$ is $l$. Observing that $z_i-z_{i+1}=\tilde{z}_i-\tilde{z}_{i+1}$, we only need to estimate $\tilde{z}_i$. From (\ref{e:vl-vi}) and $v_i-v_l=\text{curl}\,\tilde{z}_i$, we deduce that 
\begin{align*}
    \text{curl}(D_{t,l}\tilde{z}_i)=-\text{div}((\tilde{z}_i\times\nabla)(v_l+v_i^T))-\nabla(p_i-p_l)-\text{div}\,\mathring{R}_l,
\end{align*}
where we use the notation $[(z\times\nabla)v]^{ij}=\epsilon_{ikl}z^k\partial_lv^j$. Thus, we obtain the following proposition, the proof of which is exactly as \cite[Proposition 3.4]{B1}.
\begin{prop}\label{zi}
    For any $|t-t_i|\leq\tau_q$ and $N\geq0$, there holds
    \begin{align}
        \|\tilde{z}_i\|_{N+\alpha}\lesssim\tau_q\delta_{q+1}l^{-N+\alpha}\label{B(vl-vi)},
        \\
        \|z_i-z_{i+1}\|_{N+\alpha}\lesssim\tau_q\delta_{q+1}l^{-N+\alpha}.
    \end{align}
\end{prop}

\subsection{Gluing procedure}
Let
\begin{align*}
    t_i=i\tau_q,\quad I_i=[t_i+\frac{1}{3}\tau_q,t_i+\frac{2}{3}\tau_q]\cap[0,T],\quad J_i=[t_i-\frac{1}{3}\tau_q,t_i+\frac{1}{3}\tau_q]\cap[0,T].
\end{align*}
Obviously, $[0,T]$ is decomposed into pairwise disjoint intervals by $\{I_i,J_i\}$. We choose a partition of unity $\{\chi_i\}$ in time that satisfy the following three properties:
\begin{itemize}
    \item The cutoffs form a partition of unity
    \begin{gather}
    \sum_i\chi_i=1.
    \end{gather}
    \item For any i,
    \begin{align}
        \mathrm{supp}\chi_i\cap \mathrm{supp}\chi_{i+2}=\emptyset,\quad\mathrm{supp}\chi_i\subset(t_i-\frac{2}{3}\tau_q,t_i+\frac{2}{3}\tau_q),\quad\chi_i|_{J_i}=1.\label{chi}
    \end{align}
    \item For any $i$ and $N$, 
    \begin{align}
        \|\partial_t^N\chi_i\|_0\lesssim\tau_q^{-N}.
    \end{align}
\end{itemize}

We define
\begin{align*}
    \bar{v}_q=\sum_i\chi_iv_i,\quad\bar{p}_q^{(1)}=\sum_i\chi_ip_i.
\end{align*}
Obviously, $\text{div}\,\bar{v}_q=0$. Moreover, due to (\ref{chi}), we obtain
\begin{align*}
\begin{cases}
\partial_t\bar{v}_q+\text{div}(\bar{v}_q\otimes\bar{v}_q)+\nabla\bar{p}_q^{(1)}\\
=\theta_le_3+\partial_t\chi_i(v_i-v_{i+1})-\chi_i(1-\chi_i)\text{div}((v_i-v_{i+1})\otimes(v_i-v_{i+1})),\quad\forall t\in I_i,\\
\partial_t\bar{v}_q+\text{div}(\bar{v}_q\otimes\bar{v}_q)+\nabla\bar{p}_q^{(1)}=\theta_le_3,\quad\forall t\in J_i.
\end{cases}
\end{align*}

In order to construct the new Reynolds tensor $\mathring{\bar{R}}_q$, we recall the "inverse divergence" operator:
\begin{align*}
    (\mathcal{R}f)^{ij}&=\mathcal{R}^{ijk}f^k,\\
    \mathcal{R}^{ijk}&=-\frac{1}{2}\Delta^{-2}\partial_i\partial_j\partial_k-\frac{1}{2}\Delta^{-2}\partial_k\delta_{ij}+\Delta^{-1}\partial_i\delta_{jk}+\Delta^{-1}\partial_j\delta_{ik},
\end{align*}
where $f$ is a vector function with zero mean on $\mathbb{T}^3$. By direct calculation (see \cite{D1}), we have $\mathcal{R}f$ is symmetric and
\begin{align}
    \text{div}(\mathcal{R}f)=f
\end{align}
for any vector function $f$ with zero mean on $\mathbb{T}^3$.   

We define
\begin{align}
    &\mathring{\bar{R}}_q=
    \begin{cases}
        \partial_t\chi_i\mathcal{R}(v_i-v_{i+1})-\chi_i(1-\chi_i)(v_i-v_{i+1})\mathring{\otimes}(v_i-v_{i+1}),\quad\forall t\in I_i,\\
        0,\quad\forall t\in J_i.
    \end{cases}\\
    &\bar{p}_q= \bar{p}_q^{(1)}-\frac{1}{3}\chi_i(1-\chi_i)(|v_i-v_{i+1}|^2-\int_{\mathbb{T}^3}|v_i-v_{i+1}|^2).
\end{align}
Therefore, we get
\begin{align}
    \begin{cases}
\partial_t\bar{v}_q+\text{div}(\bar{v}_q\otimes\bar{v}_q)+\nabla\bar{p}_q=\theta_le_3+\text{div}\mathring{\bar{R}}_q,\\\text{div}\,\bar{v}_q=0.
    \end{cases}
\end{align}

\begin{rema}
    Here we verify $\mathcal{R}(v_i-v_{i+1})$ is well defined i.e $v_i-v_{i+1}$ has zero mean. Recall that $(v_q,\theta_q)$ solves (\ref{e:Boussinesq-Reynold}) and $\text{div}\,v_q=0$, we have 
    \begin{align*}
        \frac{d}{dt}_t\int_{\mathbb{T}^3}\theta_q=-\int_{\mathbb{T}^3}v_q\cdot\nabla\theta_q-\int_{\mathbb{T}^3}\Delta\theta_q=0.
    \end{align*}
    Since $\theta^0$ has zero mean, we obtain $\theta_q$ has zero mean for all $t$. Furthermore, $v_q$ has constant mean:
    \begin{align*}
        \frac{d}{dt}\int_{\mathbb{T}^3}v_q=-\int_{\mathbb{T}^3}({\rm div}(v_q\otimes v_q)+\nabla p_q-\theta_q e_3-{\rm div}\mathring{R_q})=0.
    \end{align*}
    Thus the mean of $v_l$ is constant. Applying the same to $v_i$, it also has constant mean. Since $v_l$ and $v_i$ coincide at the time $t_i$, we have $v_i-v_l$ has zero mean for all $t$. Due to the same reason for $v_{i+1}-v_l$ and we conclude that $v_i-v_{i+1}$ has zero mean from $v_i-v_{i+1}=(v_i-v_l)-(v_{i+1}-v_l)$.
\end{rema}

The following proposition can be easily obtained from the definition of $\bar{v}_q$, (\ref{vi-vl}) and Proposition \ref{zi}, which can be found in Proposition 4.3 and Proposition 4.4 of \cite{B1}.

\begin{prop}\label{propbarv}
    For all $N\geq0$, $\bar{v}_q$ satisfies
    \begin{gather}
        \|\bar{v}_q-v_l\|_{\alpha}\lesssim\delta_{q+1}^{\frac{1}{2}}l^{\alpha},\label{barvq-vl}\\
        \|\bar{v}_q-v_l\|_{N+\alpha}\lesssim\tau_q\delta_{q+1}l^{-N-1+\alpha},\\
        \|\bar{v}_q\|_{1+N}\lesssim\delta_q^{\frac{1}{2}}\lambda_ql^{-N}.\label{barvq}\\
        \|\mathring{\bar{R}}_q\|_{N+\alpha}\lesssim\delta_{q+1}l^{-N+\alpha}.\label{ringbarRq}
    \end{gather}
\end{prop}

In the next proposition, we will demonstrate that the energy of $\bar{v}_q$ is roughly equivalent to $v_l$.
\begin{prop}
    The difference in the energies between $\bar{v}_q$ and $v_l$ satisfies the following estimate:
    \begin{align}
        \left|\int_{\mathbb{T}^3}|\bar{v}_q|^2-|v_l|^2dx\right|\lesssim\delta_{q+1}l^{\alpha}.
    \end{align}
\end{prop}
\begin{proof}
    Since $v_l$ and $v_i$ are smooth functions satisfying (\ref{e:vl}) and (\ref{Euler with F}) respectively, we have
    \begin{align*}
        \left|\frac{d}{dt}\int_{\mathbb{T}^3}|v_i|^2-|v_l|^2dx\right|&=2\left|\int_{\mathbb{T}^3}\nabla v_l:\mathring{R}_l+\theta_l(v_i-v_l)dx\right|\\&\lesssim\|\nabla v_l\|_0\|\mathring{R}_l\|_0+\|\theta_l\|_0\|v_i-v_l\|_0\\&\lesssim\delta_q^{\frac{1}{2}}\lambda_q\delta_{q+1}+\tau_q\delta_{q+1}l^{-1+\alpha}\\&\lesssim\tau_q^{-1}\delta_{q+1}l^{\alpha}.
    \end{align*}
    where we have used (\ref{Rl}), (\ref{v1}) and (\ref{vi-vl}). Moreover, since $v_l$ and $v_i$ coincide at the time $t_i$, after integrating in time we obtain 
    \begin{align}
        \left|\int_{\mathbb{T}^3}|\bar{v}_q|^2-|v_l|^2dx\right|=\left|\int_{\mathbb{T}^3}|v_i|^2-|v_l|^2dx\right|\lesssim\delta_{q+1}l^{\alpha},\quad\forall t\in J_i.
    \end{align}
    Observe that $t\in I_i$,
    \begin{align*}
        |\bar{v}_q|^2-|v_l|^2&=|\chi_iv_i+(1-\chi_i)v_{i+1}|^2-|v_l|^2\\&=\chi_i(|v_i|^2-|v_l|^2)+(1-\chi_i)(|v_{i+1}|^2-|v_l|^2)-\chi_i(1-\chi_i)|v_i-v_{i+1}|^2,
    \end{align*}
    therefore, we deduce
    \begin{align*}
        \left|\int_{\mathbb{T}^3}|\bar{v}_q|^2-|v_l|^2dx\right|&\lesssim\left|\int_{\mathbb{T}^3}(|v_i|^2-|v_l|^2)+(|v_{i+1}|^2-|v_l|^2)dx\right|+\left|\int_{\mathbb{T}^3}|v_i-v_{i+1}|^2dx\right|\\&\lesssim\delta_{q+1}l^{\alpha}+(\tau_q\delta_{q+1}l^{-1+\alpha})^2\\&\lesssim\delta_{q+1}l^{\alpha},
    \end{align*}
    which completes the proof.
\end{proof}

\section{Perturbation}\label{Perturbation}
In this section, we will provide a brief overview of how the perturbation $w_{q+1}$ is constructed. It is based on the Mikado flows and finally $v_{q+1}$ will be defined as 
\begin{align*}
    v_{q+1}=\bar{v}_q+w_{q+1}.
\end{align*}

\subsection{Mikado flows and squiggling stripes}
We denotes $B_{\frac{1}{2}}(\text{Id})$ is the metric ball whose radius is $\frac{1}{2}$ and which is centered on Id in $\mathcal{S}_+^{3\times 3}$. We firstly recall the Mikado flows given in \cite{6}:
\begin{lemm}\label{Mikado}
    There exists a smooth vector field
    \begin{align*}
        W:B_{\frac{1}{2}}(\text{Id})\times\mathbb{T}^3\to\mathbb{R}^3
    \end{align*}
    such that for any $R\in B_{\frac{1}{2}}(\text{Id})$, it satisfies:
    \begin{gather}
        \begin{cases}
            \text{div}_\xi(W(R,\xi)\otimes W(R,\xi))=0,\\
            \text{div}_\xi\, W(R,\xi)=0,
        \end{cases}\\
        \fint_{\mathbb{T}^3}\,W(R,\xi)d\xi=0,\\
        \fint_{\mathbb{T}^3}\,W(R,\xi)\otimes W(R,\xi)d\xi=R.
     \end{gather}        
\end{lemm}

Using the fact that $\xi \to W(R,\xi)$ is $\mathbb{T}^3$-periodic with  zero mean in $\xi$, we get 
\begin{gather}
    W(R,\xi)=\sum_{k\in\mathbb{Z}^3 \backslash\{0\}}a_k(R)e^{ik\cdot\xi},\\
    W(R,\xi)\otimes W(R,\xi)=R+\sum_{k\in\mathbb{Z}^3 \backslash\{0\}}C_k(R)e^{ik\cdot\xi},
\end{gather}

where $R\to a_k(R)$ and $R\to C_k(R)$ are smooth functions satisfying $a_k(R)\cdot k=0$ and $C_k(R)k=0$. Moreover, due to the smoothness of $W$, we deduce
\begin{align}
    \sup_{R\in B_{\frac{1}{2}}(\text{Id})}|D^N_Ra_k(R)|\leq\frac{C(N,m)}{|k|^m},\quad\sup_{R\in B_{\frac{1}{2}}(\text{Id})}|D^N_RC_k(R)|\leq\frac{C(N,m)}{|k|^m}\label{Cak}
\end{align}
for any $m,N\in\mathbb{N}$.

Recalling that the support of $\mathring{\bar{R}}_q$ is contained in the set $\cup_iI_i\times\mathbb{T}^3$, we introduce the following squiggling stripes functions. More precisely, there exist smooth nonnegative cutoff functions $\{\eta_i(t,x)\}$ satisfying the following properties:
\begin{itemize}
    \item $\eta_i\in C^{\infty}([0,T]\times\mathbb{T}^3)$ with $0\leq\eta_i\leq1$.
    \item $\text{supp}\,\eta_i\cap \text{supp}\,\eta_j=\emptyset$ for $i\neq j$.
    \item $I_i\times\mathbb{T}^3\subset\{(t,x)|\eta_i=1\}$.
    \item $\text{supp}\,\eta_i\subset I_i\cup J_i\cup J_{i+1}\times\mathbb{T}^3=(t_i-\frac{1}{3}\tau_q,t_{i+1}+\frac{1}{3}\tau_q)\cap[0,T]\times\mathbb{T}^3$.
    \item There exists a positive constant $c_0>0$ such that 
    \begin{align}
        \sum_i\int_{\mathbb{T}^3}\eta_i^2(t,x)dx\geq c_0, \quad\forall t\in[0,T].\label{intrho}
    \end{align}
    \item $\|\partial_t^n\eta_i\|_m\leq C(n,m)\tau_q^{-n},\quad\forall n,m\geq0$.
\end{itemize}
The construction of $\{\eta_i(t,x)\}$ can be found in \cite[Lemma 5.3]{B1}.

\subsection{The perturbation and the constant \textit{M}}
Define 
\begin{gather*}
    \rho_q(t)\triangleq\frac{1}{3}(e(t)-\frac{\delta_{q+2}}{2}-\int_{\mathbb{T}^3}|\bar{v}_q|^2dx),\\
    \rho_{q,i}(t,x)\triangleq\frac{\eta_i^2(t,x)}{\sum_j\int_{\mathbb{T}^3}\eta_j^2(t,y)dy}\rho_q(t),\\
    R_{q,i}(t,x)\triangleq\rho_{q,i}\text{Id}-\eta_i^2\mathring{\bar{R}}_q,\\
    \widetilde{R}_{q,i}\triangleq\frac{\nabla\Phi_iR_{q,i}(\nabla\Phi_i)^T}{\rho_{q,i}},
\end{gather*}
where $\Phi_i$ is the back flows of the velocity $\bar{v}_q$ satisfying the transport equation
\begin{align*}
    \begin{cases}
        \partial_t\Phi_i+\bar{v}_q\cdot\nabla\Phi_i=0,\\
        \Phi_i(t_i,x)=x.
    \end{cases}
\end{align*}

The principal part of the perturbation $w_{q+1}$ is formulated as
\begin{align}
    w_0\triangleq\sum_i\rho_{q,i}^{\frac{1}{2}}(\nabla\Phi_i)^{-1}W(\widetilde{R}_{q,i},\lambda_{q+1}\Phi_i)=\sum_{k\neq0}\sum_i\rho_{q,i}^{\frac{1}{2}}(\nabla\Phi_i)^{-1}a_k(\widetilde{R}_{q,i})e^{i\lambda_{q+1}k\cdot\Phi_i}\nonumber,
\end{align}
where $a_k(\widetilde{R}_{q,i})$ is well defined due to (\ref{Rqi1}). In order to let $\text{div}\,w_{q+1}=0$, we define the corrector part of $w_{q+1}$ as
\begin{align}
    w_c\triangleq\frac{i}{\lambda_{q+1}}\sum_{i,k\neq0}[\text{curl}((\rho_{q,i}^{\frac{1}{2}})\frac{\nabla\Phi_i^T(k\times a_k(\widetilde{R}_{q,i}))}{|k|^2})]e^{i\lambda_{q+1}k\cdot\Phi_i}.
\end{align}

We define the perturbation
\begin{align}
    w_{q+1}\triangleq w_0+w_c=\text{curl}(\sum_{i,k\neq0}\rho_{q,i}^{\frac{1}{2}}\frac{\nabla\Phi_i^T(ik\times a_k(\widetilde{R}_{q,i}))}{\lambda_{q+1}|k|^2}e^{i\lambda_{q+1}k\cdot\Phi_i})\label{defwq+1},
\end{align}
here one can verify the second equality by direct calculation (see \cite{B1}) and thus $w_{q+1}$ is divergence-free.

\begin{lemm}\label{lemmarho}
    If $a$ is sufficiently large, for any $N\geq0$ we have
    \begin{itemize}
        \item the back flows $\{\nabla\Phi_i\}$ satisfy
        \begin{gather}
        \|\nabla\Phi_i-\text{Id}\|_0\lesssim\tau_q\delta_q^{\frac{1}{2}}\lambda_q\leq\frac{1}{10}, \quad\forall t\in \text{supp}(\eta_i),\label{梯度phi}\\
        \|\nabla\Phi_i\|_N+\|(\nabla\Phi_i)^{-1}\|_N\lesssim l^{-N},\quad\forall t\in \text{supp}(\eta_i),\label{梯度phiN}
        \end{gather}
        \item $\rho_q(t)$ satisfies
        \begin{gather}
            \frac{\delta_{q+1}}{8\lambda_q^\alpha}\leq|\rho_q(t)|\leq\delta_{q+1}\label{rho},\\
        \|\partial_t\rho_q\|_0\lesssim\delta_{q+1}\delta_q^{\frac{1}{2}}\lambda_q,
        \end{gather}
        \item $\{\rho_{q,i}\}$ satisfy
        \begin{gather}
            \|\rho_{q,i}\|_0\leq\frac{\delta_{q+1}}{c_0},\label{rhoqi的最大值估计}\\
        \|\rho_{q,i}^{\frac{1}{2}}\|_N\lesssim\delta_{q+1},\label{rho1/2}\\
        \|\partial_t\rho_{q,i}\|_N\lesssim\delta_{q+1}\tau_q^{-1}.
        \end{gather}
    \end{itemize}
    Furthermore, for all $(t,x)\in\text{supp}\,\eta_i\times\mathbb{T}^3$,
    $\widetilde{R}_{q,i}(t,x)$ is symmetric and satisfies
    \begin{gather}
        \|\widetilde{R}_{q,i}-\text{Id}\|_0\leq\frac{1}{2},\label{Rqi1}\\
        \|\widetilde{R}_{q,i}\|_N\lesssim l^{-N},\quad\forall N\geq0.\label{Rqi2}
    \end{gather}
\end{lemm}
\begin{proof}
    The proof can be found in \cite[Lemma 5.4 and Proposition 5.7]{B1}. Here we just give a proof of (\ref{rho1/2}) and (\ref{Rqi1}).
    By the definition of $\rho_{q,i}$, we have 
    \begin{align*}
        \rho_{q,i}^{\frac{1}{2}}(t,x)=\frac{\eta_i(t,x)\rho_q^\frac{1}{2}(t)}{(\sum_j\int_{\mathbb{T}^3}\eta_j^2(t,y)dy)^\frac{1}{2}},
    \end{align*}
    thus (\ref{rho1/2}) follows from (\ref{intrho}), (\ref{rho}) and the fact that $\eta_i$ is a smooth function on $[0,T]\times\mathbb{T}^3$.
    By the definition of $\widetilde{R}_{q,i}$, we have
    \begin{align*}
        \widetilde{R}_{q,i}-\text{Id}&=-\nabla\Phi_i\frac{\eta_i^2\mathring{\bar{R}}_q}{\rho_{q,i}}\nabla\Phi_i^T+\nabla\Phi_i\nabla\Phi_i^T-\text{Id}\\
        &=-\nabla\Phi_i\frac{\sum_j\int_{\mathbb{T}^3}\eta_j^2dy\,\mathring{\bar{R}}_q}{\rho_{q}}\nabla\Phi_i^T+(\nabla\Phi_i-\text{Id})\nabla\Phi_i^T+\nabla\Phi_i^T-\text{Id}.
    \end{align*}
    Applying (\ref{ringbarRq}), (\ref{梯度phi}) and (\ref{rho}), we obtain
    \begin{align*}
        \|\widetilde{R}_{q,i}-\text{Id}\|_0\leq C(l\lambda_q)^{\alpha}+\frac{11}{100}+\frac{1}{10}\leq\frac{1}{2}
    \end{align*}
    if a is sufficiently large.
\end{proof}

\begin{defi}
    The constant $M$ in Proposition \ref{Prop1} is defined as 
    \begin{align}
        M=\max\{\sqrt{\frac{M_1}{4\pi^3}},12\sum_{k\neq0}\frac{11}{9}\frac{C(0,5)}{c_0^{\frac{1}{2}}|k|^4}\}\label{M}
    \end{align}
    where $M_1$ is defined in (\ref{M1}) and $C(0,5)$ is defined in (\ref{Cak}).
\end{defi}

\begin{rema}
    Combining with (\ref{v0N}) and (\ref{M}), the starting velocity $v_0$ that we choose in Section \ref{start} satisfies (\ref{p3}) with $q=0$.
\end{rema}

\begin{prop}\label{propw}
    If $a$ is sufficiently large, there holds
    \begin{gather}
        \|w_0\|_0+\lambda_{q+1}^{-1}\|w_0\|_{1}\leq\frac{M}{4}\delta^{\frac{1}{2}}_{q+1},\label{w0}\\
        \|w_c\|_0+\lambda_{q+1}^{-1}\|w_c\|_{1}\lesssim\delta^{\frac{1}{2}}_{q+1}l^{-1}\lambda^{-1}_{q+1},\label{wc}\\
        \|w_{q+1}\|_0+\lambda_{q+1}^{-1}\|w_{q+1}\|_{1}\leq\frac{M}{2}\delta^{\frac{1}{2}}_{q+1},\label{wq+1}\\
        \|w_{q+1}\|_N\lesssim\delta^{\frac{1}{2}}_{q+1}\lambda_{q+1}^N,\quad\forall N\geq2,\label{wq+1N}
    \end{gather}
\end{prop}
\begin{proof}
    Using (\ref{梯度phi}), we conclude 
    \begin{align}
        (\nabla\Phi_i)^{-1}(t,x)\leq\frac{10}{9},\quad\nabla\Phi_i(t,x)\leq\frac{11}{10},\quad\forall (t,x)\in\text{supp}\eta_i.\label{phi的梯度的最大值估计}
    \end{align}
    Thus, by the fact that $\{\rho_{q,i}\}$ have disjoint supports, we obtain
    \begin{align}
        \|w_0\|_0\leq\sum_{k\neq0}\frac{\delta_{q+1}^{\frac{1}{2}}}{c_0^{\frac{1}{2}}}\frac{10}{9}\frac{C(B_{\frac{1}{2}}(\text{Id}),0,5)}{|k|^5}\leq\frac{M}{12}\delta_{q+1}^{\frac{1}{2}}\label{w00},
    \end{align}
    where we used (\ref{Cak}), (\ref{rhoqi的最大值估计}) and (\ref{phi的梯度的最大值估计}). Similarly, we deduce
    \begin{align*}
        \|\nabla w_0\|_0&\leq\sum_{i,k\neq0}\|\nabla(\rho_{q,i}^{\frac{1}{2}}(\nabla\Phi_i)^{-1}a_k(\widetilde{R}_{q,i}))e^{i\lambda_{q+1}k\cdot\Phi_i}\|_0+\|\rho_{q,i}^{\frac{1}{2}}(\nabla\Phi_i)^{-1}a_k(\widetilde{R}_{q,i})\lambda_{q+1}k\cdot \nabla\Phi_i e^{i\lambda_{q+1}k\cdot\Phi_i}\|_0\\
        &\leq\sum_{k\neq0}\bar{C}\delta_{q+1}^{\frac{1}{2}}l^{-1}+\frac{\delta_{q+1}^{\frac{1}{2}}}{c_0^{\frac{1}{2}}}\frac{10}{9}\frac{C(B_{\frac{1}{2}}(\text{Id}),0,5)}{|k|^4}\frac{11}{10}\\
        &\leq\sum_{k\neq0}\bar{C}\delta_{q+1}^{\frac{1}{2}}l^{-1}+\frac{M}{12}\delta_{q+1}^{\frac{1}{2}},
    \end{align*}  
    where $\bar{C}$ depends on $\beta,\alpha$ and $M$ but not on $a$. By the definition of $l$ in (\ref{l}), we get
    \begin{align*}
        (l\lambda_{q+1})^{-1}=\lambda_q^{1-\beta+\frac{3\alpha}{2}-b(1-\beta)}.
    \end{align*}
    Choosing $\alpha$ is sufficiently small such that $b>\frac{1-\beta+\frac{3\alpha}{2}}{1-\beta}$ and $a$ is sufficiently large, we achieve
    \begin{align}
        \|\nabla w_0\|_0\leq\frac{M}{8}\delta_{q+1}^{\frac{1}{2}}.\label{w01}
    \end{align}
    Combining (\ref{w00}) and (\ref{w01}), we conclude (\ref{w0}). Furthermore, (\ref{wc}) can be proved similarly and (\ref{wq+1}) follows as a direct consequence of (\ref{w0}) and (\ref{wc}). 
    By the definition of $w_{q+1}$ in (\ref{defwq+1}) and Lemma \ref{lemmarho}, for any $N\geq2$ we have
    \begin{align*}
        \|w_{q+1}\|_N&\lesssim\sum_{i,k\neq0}\|\rho_{q,i}^{\frac{1}{2}}\frac{\nabla\Phi_i^T(ik\times a_k(\widetilde{R}_{q,i}))}{\lambda_{q+1}|k|^2}e^{i\lambda_{q+1}k\cdot\Phi_i}\|_{N+1}\\
        &\lesssim\lambda_{q+1}^{-1}\sum_{i,k\neq0}(\|\rho_{q,i}^{\frac{1}{2}}\frac{\nabla\Phi_i^T(ik\times a_k(\widetilde{R}_{q,i}))}{|k|^2}\|_{N+1}+\|\rho_{q,i}^{\frac{1}{2}}\frac{\nabla\Phi_i^T(ik\times a_k(\widetilde{R}_{q,i}))}{|k|^2}\|_0\|e^{i\lambda_{q+1}k\cdot\Phi_i}\|_{N+1})\\
        &\lesssim\lambda_{q+1}^{-1}\delta_{q+1}^{\frac{1}{2}}(l^{-N-1}+\lambda_{q+1}^{N+1})\\
        &\lesssim\delta_{q+1}^{\frac{1}{2}}\lambda_{q+1}^{N},
    \end{align*}
    which completes the proof.
\end{proof}

\begin{prop}
    Assuming $a$ is sufficiently large, $v_{q+1}$ satisfies the following estimates:
    \begin{gather}
        \|v_{q+1}\|_0\leq C-\delta_{q+1}^{\frac{1}{2}},\\
        \|v_{q+1}\|_1\leq M\delta^{\frac{1}{2}}_{q+1}\lambda_{q+1},\\
        \|v_{q+1}\|_N\leq C'(N)\delta^{\frac{1}{2}}_{q+1}\lambda_{q+1}^N,\quad\forall N\geq2,\\
         \|v_{q+1}-v_{q}\|_0+\lambda_{q+1}^{-1}\|v_{q+1}-v_{q}\|_1\leq M\delta_{q+1}^{\frac{1}{2}}\label{零次和一次的差},
    \end{gather}
    where $\{C'(N)\}_{N\geq2}$ depend on $\beta,\alpha$ and $M$ but not on $a$ and $q$.
\end{prop}
\begin{proof}
    Combining Proposition \ref{propvl-vq}, Proposition \ref{propbarv} and Proposition \ref{propw}, we can easily get the result.
\end{proof}

\begin{defi}
    The series $\{C(N)\}_{N\geq2}$ in Proposition \ref{Prop1} is defined as \begin{align}
    C(N)=max\{C'(N),\sqrt{\frac{M_1}{4\pi^3}}\}.
\end{align}
\end{defi}

Obviously, the starting velocity $v_0$ and $v_{q+1}$ both satisfy (\ref{p4}). Moreover, as in \cite[Proposition 6.2]{B1}, we have the following energy estimate: 
\begin{prop}
    Assuming $a$ is sufficiently large, we have
    \begin{align*}
        \delta_{q+2}\lambda_{q+1}^{-\alpha}\leq e(t)-\int_{{\mathbb{T}}^3}|v_{q+1}(t,x)|^2dx\leq \delta_{q+2},\quad\forall t\in[0,T]
    \end{align*}
\end{prop}

\section{Construction of the new temperature}\label{温度构造}
In this section, we focus on the construction of a new temperature $\theta_{q+1}$. In order to obtain the estimates of $\theta_q$ in the Sobolev space, we make an additional assumption on $v_q$, that is, $v_q$ satisfies (\ref{p4}). Combining (\ref{p3}) with (\ref{p4}), we have
\begin{align}
    \|v_q\|_N\leq C(N)\delta_q^{\frac{1}{2}}\lambda_q^N, \quad\forall N\geq1,\label{vqN}
\end{align}
where $C(1)=M$. By the construction on $v_{q+1}$, we know $v_{q+1}$ also satisfies the above estimate with $q$ replaced by $q+1$.

Let $\theta_{q+1}$ be the smooth solution of the following transport-diffusion equation:
    \begin{align}
        \begin{cases}
            \partial_t\theta_{q+1}+v_{q+1}\cdot\nabla\theta_{q+1}-\Delta\theta_{q+1}=0,\\
            \theta(0,x)=\theta^0(x_3).
        \end{cases}
    \end{align}
Furthermore, the energy equality (\ref{p8}) with $q$ replaced by $q+1$ can be easily verified. 
\begin{prop}
    For any $N\geq2$, we have
    \begin{gather}
        \|\nabla\theta_{q}\|_{L^\infty L^2}\lesssim1,\label{thetaq1}\\
         \|\nabla\theta_{q+1}\|_{L^\infty L^2}\lesssim1,\label{theaq+11}\\
        \|\nabla^N\theta_{q}\|_{L^\infty L^2}\lesssim\delta_q^{\frac{1}{2}}\lambda_q^{N-1},\\
        \|\nabla^N\theta_{q+1}\|_{L^\infty L^2}\lesssim\delta_{q+1}^{\frac{1}{2}}\lambda_{q+1}^{N-1}.\label{thetaq+1N}
    \end{gather}
\end{prop}
\begin{proof}
    Since $\theta_{q}$ is a solution of $(\ref{e:Boussinesq-Reynold})_3$, we obtain
    \begin{align*}
        \|\nabla^N\theta_{q}\|_{L^\infty L^2}&\leq\|\nabla^N\theta_0\|_{L^\infty L^2}+C(N,\|v_q\|_0,\|\theta_0\|_{L^2})\sum_{k=1}^{N-1}\|\nabla ^k v_q\|_0^{\frac{N}{k+1}}\\&\lesssim\sum_{k=1}^{N-1}(\delta^{\frac{1}{2}}_{q+1}\lambda_{q+1}^k)^{\frac{N}{k+1}}\\&\lesssim\delta^{\frac{1}{2}}_{q+1}\lambda_{q+1}^{N-1}, 
    \end{align*}
    where we apply (\ref{vqN}) and (\ref{thetaN}). Furthermore, (\ref{thetaq1}) follows as a direct consequence of (\ref{p2}) and (\ref{theta1}). Arguing in a similar way, we can also obtain (\ref{theaq+11}) and (\ref{thetaq+1N}) for $\theta_{q+1}$.
\end{proof}

\begin{prop}
    \begin{align}
         \|(\theta_{q+1}-\theta_{q})(t,\cdot)\|^2_{L^2}+\int_0^t\|\nabla(\theta_{q+1}-\theta_{q})(s,\cdot)\|^2_{L^2}ds\leq C\delta_{q+1}^{\frac{1}{2}},\quad\forall t\in[0,T].
    \end{align}
\end{prop}
\begin{proof}
    By direct calculation, we obtain
    \begin{align*}
        \begin{cases}
        \partial_t(\theta_{q+1}-\theta_{q})+v_{q+1}\cdot\nabla(\theta_{q+1}-\theta_{q})-\Delta(\theta_{q+1}-\theta_{q})=(v_q-v_{q+1})\cdot\nabla\theta_q,\\
        \theta_{q+1}-\theta_{q}|_{t=0}=0.
        \end{cases}
    \end{align*}
    Using (\ref{零次和一次的差}) and (\ref{thetaq1}), a direct energy estimate give
    \begin{align}
        \|(\theta_{q+1}-\theta_{q})(t,\cdot)\|^2_{L^2}+\int_0^t\|\nabla(\theta_{q+1}-\theta_{q})(s,\cdot)\|^2_{L^2}ds\lesssim\|((v_q-v_{q+1})\cdot\nabla\theta_q\|_{L^2}\leq C\delta_{q+1}^{\frac{1}{2}},
    \end{align}
    which completes the proof.
\end{proof}

Thus, $\theta_{q+1}-\theta_{q}$ satisfies (\ref{p7}), but to estimate the new Reynolds stress $\mathring{R}_{q+1}$ that will be defined in the next section, we need to get more precise estimates of $\theta_{q+1}-\theta_{q}$.

\begin{prop}\label{theta的0.5阶导}
    \begin{align}
        \|\theta_{q+1}-\theta_{q}\|_{L^2}\lesssim(\delta_q^{\frac{1}{2}}\lambda_q)^\alpha l^{1-\alpha},\label{总的差}\\
        \|\nabla\theta_{q+1}-\nabla\theta_{q}\|_{L^2}\lesssim\delta_{q+1}^{\frac{1}{2}}.\label{导数的差}
    \end{align}
\end{prop}
\begin{proof}
    We define $\{f_1,f_2,f_3\}$ which satisfy the following equation with zero initial datum:
    \begin{align}
         \partial_tf_1+v_{q+1}\cdot\nabla f_1-\Delta f_1=(v_q-v_l)\cdot\nabla\theta_q,\label{f1}\\
         \partial_tf_2+v_{q+1}\cdot\nabla f_2-\Delta f_2=(v_l-\bar{v}_q)\cdot\nabla\theta_q,\\
         \partial_tf_3+v_{q+1}\cdot\nabla f_3-\Delta f_3=-w_{q+1}\cdot\nabla\theta_q.
    \end{align}
    Obviously we have $\theta_{q+1}-\theta_{q}=f_1+f_2+f_3$.

{\bf Step 1:}
    Since $v_q-v_l$ has zero mean, we have
    \begin{align*}
        v_q-v_l=\text{curl}\,\mathcal{B}(v_q-v_l)
    \end{align*}
    where $\mathcal{B}$ is the Biot-Savart operator. Moreover, we can write 
    \begin{align*}
        (v_q-v_l)\cdot\nabla\theta_q=(\text{curl}\,\mathcal{B}(v_q-v_l))\cdot\nabla\theta_q=\text{div}((\mathcal{B}(v_q-v_l)\times\nabla)\theta_q),
    \end{align*}
    here we use notation $[(z\times\nabla)\theta]^i=\epsilon_{ikl}z^k\partial_l\theta$ for vector fields $z$ and scalar function $\theta$.
    
    Taking $L^2$ inner products with $f_1$ in (\ref{f1}) and integrating by parts, we deduce that
    \begin{align}
        \frac{1}{2}\frac{d}{dt}\|f_1(t,\cdot)\|^2_{L^2}+\|\nabla f_1(t,\cdot)\|^2_{L^2}&=\left|\int_{\mathbb{T}^3}(\mathcal{B}(v_q-v_l)\times\nabla)\theta_q\cdot\nabla f_1dx\right|\label{f1的不等式}\\
        &\leq 4\|(\mathcal{B}(v_q-v_l)\times\nabla)\theta_q\|_{L^2}^2+\frac{1}{4}\|\nabla f_1(t,\cdot)\|_{L^2}^2.\nonumber
    \end{align}
    Additionally, using the fact that $\nabla\mathcal{B}$ is a bounded operator on H\"{o}lder space, we have
    \begin{align}
        \|\mathcal{B}(v_q-v_l)\|_0\lesssim\|\mathcal{B}v_q-\mathcal{B}v_q*\phi_l\|_0\lesssim\|\nabla\mathcal{B}v_q\|_\alpha\, l^{1-\alpha}
       \lesssim\|v_q\|_\alpha l^{1-\alpha}\lesssim(\delta_q^{\frac{1}{2}}\lambda_q)^\alpha l^{1-\alpha}.\label{B的估计}
    \end{align}
    Thus, using (\ref{thetaq1}), (\ref{f1的不等式}) and (\ref{B的估计}) we obtain
    \begin{align}
       \|f_1(t,\cdot)\|_{L^2}\lesssim\|(\mathcal{B}(v_q-v_l)\times\nabla)\theta_q\|_{L^2}\lesssim(\delta_q^{\frac{1}{2}}\lambda_q)^\alpha l^{1-\alpha}\label{f1的差}.
    \end{align}
    Applying $\nabla$ on the both side of (\ref{f1}) and a direct energy estimate give
    \begin{align}
        \|\nabla f_1(t,\cdot)\|_{L^2}\lesssim\|(v_q-v_l)\cdot\nabla\theta_q\|_{L^2}\lesssim\delta_{q+1}^{\frac{1}{2}}\lambda_q^{-\alpha}\label{f1导数的差}.
    \end{align}

{\bf Step 2:}
    By the definition of $\bar{v}_q$, we have
    \begin{align*}
        v_l-\bar{v}_q=\sum_i\chi_i( v_l-v_i)=\sum_i\chi_i\text{curl}\,\mathcal{B}( v_l-v_i).
    \end{align*}
    Considering $\{\chi_i\}$ is a partition of unity which has almost disjoint supports, arguing in a similar way as Step 1, we obtain
    \begin{gather}
        \|f_2(t,\cdot)\|_{L^2}\lesssim\|(\mathcal{B}(v_l-v_i)\times\nabla)\theta_q\|_{L^2}\lesssim\tau_q\delta_{q+1}l^{\alpha},\label{f2的差}\\\|\nabla f_2(t,\cdot)\|_{L^2}\lesssim\|(v_l-\bar{v}_q)\cdot\nabla\theta_q\|_{L^2}\lesssim\delta_{q+1}^{\frac{1}{2}}l^\alpha\label{f2导数的差},
    \end{gather}
    where we used (\ref{B(vl-vi)}) and (\ref{barvq-vl}).

{\bf Step 3:}  
    We define
     \begin{align}
         d_{i,k}=&\rho_{q,i}^{\frac{1}{2}}(\nabla\Phi_i)^{-1}a_k(\widetilde{R}_{q,i})\nabla\theta_qe^{i\lambda_{q+1}k\cdot(\Phi_i-x)}\\ \notag &+\frac{i}{\lambda_{q+1}}[\text{curl}((\rho_{q,i}^{\frac{1}{2}})\frac{\nabla\Phi_i^T(k\times a_k(\widetilde{R}_{q,i}))}{|k|^2})]\nabla\theta_qe^{i\lambda_{q+1}k\cdot(\Phi_i-x)},\label{dik的定义}
     \end{align}
     thus by the definition of $w_{q+1}$, we have 
     \begin{align}      w_{q+1}\cdot\nabla\theta_q=\sum_{k\neq0}\sum_id_{i,k}e^{i\lambda_{q+1}k\cdot x}.
     \end{align}
     Applying (\ref{复合}), (\ref{梯度phi}) and (\ref{梯度phiN}), for any $N$ sufficiently large, we deduce
     \begin{align}
         \|e^{i\lambda_{q+1}k\cdot(\Phi_i-x)}\|_N&\lesssim\lambda_{q+1}|k|l^{-N+1}+(\lambda_{q+1}|k|\tau_q\delta_q^{\frac{1}{2}}\lambda_q)^N\label{指数的估计}\\
         &\lesssim\lambda_{q+1}^{N-1}|k|^N\nonumber.
     \end{align}
     Then using Lemma \ref{lemmarho}, (\ref{dik的定义}) and (\ref{指数的估计}), we obtain
     \begin{align*}
         \|\nabla^Nd_{i,k}\|_{L^2}&\lesssim(\delta_{q+1}^{\frac{1}{2}}(l^{-N}+\delta_q^{\frac{1}{2}}\lambda_q^N+\lambda_{q+1}^{N-1})+\lambda_{q+1}^{-1}\delta_{q+1}^{\frac{1}{2}}l^{-1}(l^{-N}+\lambda_{q+1}^{-1}\delta_q^{\frac{1}{2}}\lambda_q^N+\lambda_{q+1}^{N-1}))|k|^{-m}\\
         &\lesssim\delta_{q+1}^{\frac{1}{2}}\lambda_{q+1}^{N-1}|k|^{-m},
     \end{align*}
     for any $N$ sufficiently large and any $m>0$.
     Therefore, using Lemma \ref{c2}:
     \begin{align}
    \|f_3(t,\cdot)\|_{L^2}&\lesssim\sum_{k\neq0}\frac{\|d_{i,k}\|_{L^2}}{\lambda_{q+1}|k|}+\frac{\|\nabla^Nd_{i,k}\|_{L^2}}{(\lambda_{q+1}|k|)^N}\label{f3的差}\\
          &\lesssim\delta_{q+1}^{\frac{1}{2}}\lambda_{q+1}^{-1}\nonumber,
     \end{align}
     here we fix a sufficiently large $N$. Furthermore, we have
     \begin{align}
         \|\nabla f_3(t,\cdot)\|_{L^2}\lesssim\|w_{q+1}\cdot\nabla\theta_q\|_{L^2}\lesssim\delta_{q+1}^{\frac{1}{2}}.\label{f3导数的差}
     \end{align}
     Finally, the estimate (\ref{总的差}) follows as a consequence of (\ref{f1的差}), (\ref{f2的差}) and (\ref{f3的差}), the estimate (\ref{导数的差}) follows as a consequence of (\ref{f1导数的差}), (\ref{f2导数的差}) and (\ref{f3导数的差}).
\end{proof}

\section{Construction of the new Reynolds stress error}
The new Reynolds stress is defined by
\begin{align*}
    \mathring{R}_{q+1}\triangleq I_1+I_2+I_3
\end{align*}
where 
\begin{gather}
    I_1\triangleq\mathcal{R}(w_{q+1}\cdot\nabla\bar{v}_q)+\mathcal{R}(\partial_tw_{q+1}+\bar{v}_q\cdot\nabla w_{q+1})+\mathcal{R}\text{div}(-\sum_iR_{q,i}+w_{q+1}\otimes w_{q+1}),\\
     I_2\triangleq\mathcal{R}((\theta_q-\theta_{q+1})e_3), \quad I_3\triangleq\mathcal{R}((\theta_l-\theta_q)e_3).
\end{gather}
The new pressure is defined by
\begin{align*}
    p_{q+1}(t,x)\triangleq\bar{p}_q(t,x)-\sum_i\rho_{q,i}(t,x)+\rho_q(t).
\end{align*}
With the above definition and the definition of $v_{q+1}$, we have
\begin{align*}
    \begin{cases}
    \partial_tv_{q+1}+\text{div}(v_{q+1}\otimes v_{q+1})+\nabla p_{q+1}=\theta_{q+1} e_3+\text{div}\mathring{R}_{q+1}, \quad\quad \\
    \text{div}\,v_q=0.
    \end{cases}
\end{align*}

As in \cite[Proposition 6.1]{B1}, we obtain 
\begin{align}
    \|I_1\|_0\lesssim\frac{\delta_{q+1}^{\frac{1}{2}}\delta_{q}^{\frac{1}{2}}\lambda_q}{\lambda_{q+1}^{1-4\alpha}}\leq\frac{1}{3}\delta_{q+2}\lambda_{q+1}^{-3\alpha},\label{I1}
\end{align}
here we omit the proof. In the following, we focus on the estimate of $I_2$ and $I_3$.

\begin{lemm}\label{lemma of R}
    For any $s>\frac{1}{2}$, we have
    \begin{align}
        \|\mathcal{R}(v)\|_0\leq C_s\|v\|_{\dot{H}^s}.
    \end{align}
\end{lemm}
\begin{proof}
    Let $v(x)=\sum_{k\in\mathbb{Z}^3}v_ke^{ik\cdot x}$ for any $x\in\mathbb{T}^3$. By the definition of $\mathcal{R}$, it holds that
    \begin{align*}
        \mathcal{R}(v)(x)=\sum_{k\in\mathbb{Z}^3,k\neq0}(\frac{-iv_k\otimes k}{|k|^2}+\frac{-ik\otimes v_k}{|k|^2}+\frac{iv_k\cdot k}{|k|^2})e^{ik\cdot x}.
    \end{align*}
Thus for any $s>\frac{1}{2}$, we have
\begin{align*}
    \|\mathcal{R}(v)\|_0\leq C\sum_{k\in\mathbb{Z}^3,k\neq0}\frac{|v_k|}{|k|}\leq C(\sum_{k\in\mathbb{Z}^3,k\neq0}\frac{1}{|k|^{2(1+s)}})^{\frac{1}{2}}(\sum_{k\in\mathbb{Z}^3,k\neq0}|v_k|^2|k|^{2s})^{\frac{1}{2}}\leq C_s\|v\|_{\dot{H}_s}.
\end{align*}
\end{proof}

\begin{prop}
    Assuming $a$ is sufficiently large, we have
    \begin{align}
        \|I_2\|_0\leq\frac{1}{3}\delta_{q+2}\lambda_{q+1}^{-3\alpha},\label{I2}\\\|I_3\|_0\leq\frac{1}{3}\delta_{q+2}\lambda_{q+1}^{-3\alpha}.\label{I3}
    \end{align}
\end{prop}
\begin{proof}
    Let $s>\frac{1}{2}$, using Lemma \ref{lemma of R}, (\ref{插值}) and Proposition \ref{theta的0.5阶导}, we obtain
    \begin{align*}
        \|I_2\|_0\leq\|\mathcal{R}(\theta_{q+1}-\theta_q)\|_0\lesssim C_s\|\theta_{q+1}-\theta_q\|_{\dot{H}^s}\lesssim((\delta_q^{\frac{1}{2}}\lambda_q)^\alpha l^{1-\alpha})^{1-s}\delta_{q+1}^{\frac{1}{2}s}\leq\frac{1}{3}\delta_{q+2}\lambda_{q+1}^{-3\alpha}.
    \end{align*}
    To obtain the parameter inequality
    \begin{align}
        ((\delta_q^{\frac{1}{2}}\lambda_q)^\alpha l^{1-\alpha})^{1-s}\delta_{q+1}^{\frac{1}{2}s}\leq\frac{1}{3}\delta_{q+2}\lambda_{q+1}^{-3\alpha},\label{参数不等式}
    \end{align}
    we divides by the right-side, take logarithms, divides by $\text{log}\,\lambda_q$  and let $\alpha$ tend to zero, we have to ensure
    \begin{align*}
        -b\beta s+(-b\beta+\beta-1)(1-s)+2b^2\beta<0.
    \end{align*}
    Let $s$ tend to $\frac{1}{2}$, the right side become
    \begin{align*}
        2\beta b^2-b\beta+\frac{1}{2}\beta-1.
    \end{align*}
    A direct calculation yields
    \begin{align*}
        2\beta b^2-b\beta+\frac{1}{2}\beta-1<0, \quad\forall 1<b<\frac{\beta+\sqrt{4\beta-3\beta^2}}{4\beta}.
    \end{align*}
    Thus (\ref{参数不等式}) holds when $\alpha$ is sufficiently small,  $s$ is close to $\frac{1}{2}$ and a is sufficiently large.
    
    By the fact that $\mathcal{R}$ is a operator of degree $-1$, we obtain
\begin{align}
    \|I_3\|_0\lesssim\|\mathcal{R}\theta_q*\phi_l-\mathcal{R}\theta_q\|_0\lesssim\|\mathcal{R}\theta_q\|_{1-\alpha}l^{1-\alpha}\lesssim\|\theta_q\|_0l^{1-\alpha}\lesssim l^{1-\alpha}\leq\frac{1}{3}\delta_{q+2}\lambda_{q+1}^{-3\alpha},
\end{align}
here the proof of the last parameter inequality is similar to the above.
\end{proof}

Finally, by (\ref{I1}), (\ref{I2}) and (\ref{I3}), the new stress error $\mathring{R}_{q+1}$ satisfies (\ref{p1}) with $q$ replaced with $q+1$:
\begin{align*}
    \|\mathring{R}_{q+1}\|_0\leq\delta_{q+2}\lambda_{q+1}^{-3\alpha}.
\end{align*}

\appendix

\section{H\"{o}lder space and Sobolev space}\label{appendix A}
In the following, $m=0,1,...,\alpha\in(0,1)$ and $\gamma$ is a multi-index. Firstly, we denote
\begin{align*}
    \|f\|_0:=\sup_{t,x}|f(t,x)|.
\end{align*}
The usual spatial H\"{o}lder seminorms are introduced as follows.
\begin{gather*}
   [f]_m=\max_{|\gamma|=m}\|D^\gamma f\|_0,\\
   [f]_{m+\alpha}=\max_{|\gamma|=m}\sup_{x \neq y,t}\frac{|D^\gamma f(t,x)-D^\gamma f(t,y)|}{|x-y|^\alpha} ,
\end{gather*}    
where $D^\gamma$ are space derivatives only. We define the H\"{o}lder norms as
\begin{align}
   \|f\|_{m+\alpha}=\sum^{m}_{k=0}[f]_k,\|f\|_{m+\alpha}=\|f\|_{m}+[f]_{m+\alpha}\nonumber.
\end{align}     

Next, we recall the following classical inequalities:
\begin{gather}
    [fg]_N\leq C([f]_N\|g\|_0+[g]_N\|f\|_0), \quad\forall N\in \mathbb{N},\\
    \|fg\|_N\leq C(\|f\|_N\|g\|_0+\|g\|_N\|f\|_0),\quad\forall N\in \mathbb{N},\label{乘法分配}\\
    [f\circ g]_N\leq C([f]_1[g]_N+\|\nabla f\|_{N-1}[g]_1^N),\quad\forall N\in \mathbb{N}_+.\label{复合}
\end{gather}

we also recall the quadratic commutator estimate (see \cite{cons}):
\begin{align}
    \|(f*\phi_l)(g*\phi_l)-fg*(\phi_l)\|_N\leq Cl^{2-N}\|f\|_1\|g\|_1,\quad\forall N\geq0,\label{commu} 
\end{align}
where $\phi$ is a standard mollifier and $\phi_l=l^{-3}\phi(l^{-1}\cdot)$. 

We define the norm of homogeneous Sobolev space:
\begin{align*}
    \|v\|_{\dot{H}^s}^2=\sum_{k\in\mathbb{Z}^3,k\neq0}|v_k|^2|k|^{2s},
\end{align*}
where 
\begin{align*}
    v(x)=\sum_{k\in\mathbb{Z}^3}v_ke^{ik\cdot x},\quad\forall x\in\mathbb{T}^3.
\end{align*}

Recall the following interpolation inequality:
\begin{align}
    \|v\|_{\dot{H}^s}\lesssim \|v\|_{\dot{H}^{s_1}}^{\alpha}\|v\|_{\dot{H}^{s_2}}^{1-\alpha}\label{插值},
\end{align}
where $s=\alpha s_1+(1-\alpha) s_2$.

\section{Estimates for the transport equation}
In this appendix, we recall some classical estimates for the smooth solution of the transport equation:
\begin{align}
    \begin{cases}\label{e:trans}
        \partial_tf+v\cdot\nabla f=g,\\
        f|_{t=0}=f_0,\\
    \end{cases}
\end{align}
where $v(t,x)$ is a given smooth vector field and $g$ is also smooth.
\begin{prop}
    Let $f$ be the solution of (\ref{e:trans}), then $f$ satisfies
    \begin{gather}
        \|f(t)\|_0\leq\|f_0\|_0+\int_0^t\|g(s,\cdot)\|_0ds,\\
        \|f(t)\|_{\alpha}\leq (\|f_0\|_{\alpha}+\int_0^t\|g(s,\cdot)\|_{\alpha}ds)\text{exp}(\int_0^t\|\nabla v(s,\cdot)\|_0ds), \quad\forall \alpha\in(0,1).\label{trans1}
     \end{gather}
\end{prop}

\section{Estimates for the transport-diffusion equation}
In this section, we recall the following energy inequality for the transport-diffusion equation, and the proof can be found in \cite[Lemma 3.7]{L1}.
\begin{lemm}
    Let $\theta_0(x)$ be a smooth function on $\mathbb{T}^3$ and $v(t,x)$ be a smooth velocity filed. Suppose $\theta$ is a smooth solution of the following transport-diffusion equation:
    \begin{align}
        \begin{cases}
            \partial_t\theta+v\cdot\nabla\theta-\Delta\theta=0,\\
            \text{div}\,v=0,\\
            \theta(0,x)=\theta_0(x).
        \end{cases}
    \end{align}
    Then it holds that 
    \begin{gather}
        \|\nabla\theta\|_{L^\infty L^2}\leq \|\nabla\theta_0\|_{L^\infty L^2}+C\|v\|_0^2\|\theta_0\|_{L^2},\label{theta1}\\
        \|\nabla^N\theta\|_{L^\infty L^2}\leq \|\nabla^N\theta_0\|_{L^\infty L^2}+C(N,\|v\|_0,\|\theta_0\|_{L^2})\sum_{k=1}^{N-1}\|\nabla ^kv\|_0^{\frac{N}{k+1}},\quad\forall N\geq2.\label{thetaN}
    \end{gather}
\end{lemm}

\begin{lemm}\label{c2}
    Let $v(t,x)$ be a smooth velocity field and $g(t,x)$ be a smooth function. Suppose $\theta$ is a smooth solution of the following transport-diffusion equation:
    \begin{align}
        \begin{cases}
            \partial_t\theta+v\cdot\nabla\theta-\Delta\theta=g(t,x)e^{i\lambda k\cdot x},\\
            \text{div}\,v=0,\\
            \theta(0,x)=0,
        \end{cases}
    \end{align}
    where $k$ is a vector satisfying $|k|=1$.
    Then it holds that
    \begin{align}
        \|\theta\|_{L^\infty L^2}\lesssim\frac{\|g\|_{L^\infty L^2}}{\lambda}+\frac{\|\nabla^Ng\|_{L^\infty L^2}}{\lambda^N},\quad\forall N\geq2.
    \end{align}
\end{lemm}

\smallskip
\noindent\textbf{Acknowledgments} This work was partially supported by the National Natural Science Foundation of China (No.12171493). 


\smallskip
\noindent\textbf{Data availability} The data that support the findings of this study are available on citation. The data that support the findings of this study are also available from the corresponding author upon reasonable request.

\phantomsection
\addcontentsline{toc}{section}{\refname}
\bibliographystyle{abbrv} 
\bibliography{1}

\end{document}